\documentclass{article}
\usepackage{amssymb}
\usepackage{amsmath}
\usepackage{enumerate}
\usepackage{amscd}
\usepackage{graphicx}
\usepackage{latexsym}
\usepackage{epsfig}
\usepackage[latin1]{inputenc}

\def\div{\mathfrak{Div}}
\def\deg{\mathfrak{Deg}}

\def\r{\mathbb{R}}
\def\n{\mathbb{N}}
\def\c{\mathbb{C}}
\def\q{\mathbb{Q}}
\def\s{\mathbb{S}}
\def\d{\mathbb{D}}

\def\z{\mathbb{Z}}

\def\pg{\mathfrak{p}}

\def\Fg{\mathfrak{F}}

\def\Wg{\mathfrak{W}}
\def\qg{\mathfrak{q}}

\newenvironment{proof}{\trivlist
\item[\hskip\labelsep{\em Proof}\,:]}{\hfill{$\Box$}\endtrivlist}

\title{\huge Exotic Minimal Surfaces}
\author{\Large Francisco J. L\'{o}pez   \thanks{Research partially
supported by MCYT-FEDER research project MTM2007-61775 and Junta de Andalucia Grant P06-FQM-01642.
Complete Minimal Surfaces, Parabolic Riemann Surfaces, Runge's Theorem, Connected Sum Constructions.}}

\newtheorem{lemma}{Lemma}[section]
\newtheorem{remark}{Remark}[section]
\newtheorem{theorem}{Theorem}[section]
\newtheorem{corollary}{Corollary}[section]
\newtheorem{assertion}{Claim}[section]
\newtheorem{definition}{Definition}[section]
\begin{document}
\maketitle

\begin{abstract}

We prove a general fusion theorem for complete orientable minimal surfaces in $\r^3$ with finite total curvature. As a consequence, complete orientable minimal surfaces of weak finite total curvature with exotic geometry are produced. More specifically, universal surfaces (i.e., surfaces from which all minimal surfaces can be recovered) and  space-filling surfaces with arbitrary genus and no symmetries. 
\end{abstract}
\section{Introduction} \label{sec:intro}

 As usual, a surface is said to be {\em open} if it is non-compact and has empty  boundary. An  open Riemann surface is said to be {\em hyperbolic} if an only if it carries a negative non-constant subharmonic function. Otherwise, it is said to be {\em parabolic}. Compact Riemann surfaces with empty boundary are said to be {\em elliptic}.

Let $M$ be  a Riemann surface  possibly with non empty  compact boundary, and let $X:M \to \r^3$ be a conformal minimal immersion. Throughout this paper we will always assume that $X$ extends as a conformal minimal immersion  to an open Riemann surface containing $M$ as a proper subset. When $X$ has finite total curvature (FTC for short), Huber and Osserman theorems \cite{huber,osserman} imply that $M$  has finite conformal type and the Weierstrass data of $X$ extend meromorphically to the topological ends of $M.$ This simply means that $M=M^c-E,$ where $M^c$ is a compact Riemann surface and $E\subset M^c-\partial(M^c)$ is a  finite set, and the Weierstrass data of $X$ have no essential singularities at points of $E.$ The Riemann surface $M^c$ is called the {\em Osserman compactification} of $M.$ Likewise, a conformal complete minimal immersion $X:M \to \r^3$ is said to be of {\em weak finite total curvature} (WFTC for short) if $X|_{\Omega}$ has finite total curvature (FTC for short) for any proper region $\Omega \subset M$ with compact boundary and finite topology. 

Unlike the FTC case, there exists orientable complete minimal surfaces of WFTC with arbitrary topology and conformal type (see \cite{approx}). The aim of this paper is to present some examples of this kind of surfaces with {\em exotic geometry}. 

An interesting question is whether there exists a  complete minimal surface from which all minimal surfaces could be recovered. Given an open Riemann surface $N,$ a complete conformal minimal immersion  $Y:N \to \r^3$  is said to be {\em universal} if it passes by all compact minimal surfaces  in $\r^3.$ In other words, if {\em for any} compact Riemann surface $M$ with  $\partial(M)\neq \emptyset$  and {\em any} conformal minimal immersion  $X:M\to\r^3,$  there is a sequence $\{M_n\}_{n \in \n}$ of regions in $N$ and biholomorphisms $h_n:M \to M_n,$ $n \in \n,$ such that $\{Y \circ h_n\}_{n \in \n} \to X$ uniformly on $M.$ Our first result provides an affirmative answer to the this  question (see Theorem \ref{th:universal}):
\begin{quote} 
{\bf Theorem I} {\em There exist parabolic universal minimal surfaces of WFTC.} 
\end{quote}
Any universal minimal immersion $Y:N \to \r^3$ is  {\em space-filling} (that is to say, $\overline{Y(N)}=\r^3$). As far as the author knows,  all previously known space-filling complete minimal surfaces are hyperbolic, simply connected and highly symmetric. The reason why is that their construction is based either in Schwarzian reflection on a fundamental compact domain or in the classical Björling problem (see \cite{mira} for a good setting). However, in Corollary \ref{co:space-filling}  we have shown that:
\begin{quote} {\bf Theorem II} {\em There exists space-filling complete minimal surfaces with WFTC, arbitrary (possibly infinity) genus, parabolic conformal type and no symmetries.}
\end{quote}

Both above results are based on a general connected sum construction (or fusion theorem) for complete minimal surfaces with FTC. For a thorough exposition of the details, the following notations are required.

Given two Riemann surfaces $M$ and $M^*$ possibly with non empty boundary,  $M^*$ is said to be an  {\em extension} of $M$ if $M$ is a proper subset of $M^*,$ $M\cap \partial(M^*)=\emptyset$ and $M^*-M^\circ$ contains no compact connected components that are disjoint from $\partial(M^*),$ where $M^\circ$ is the topological interior of $M$ in $M^*.$    If $M^*$ is an extension of $M$ and $j:M \to M^*$ is the inclusion, then $j_*:{\cal H}_1(M,\r) \to {\cal H}_1(M^*,\r)$ is a group monomorphism, hence up to natural identifications ${\cal H}_1(M,\r) \subset {\cal H}_1(M^*,\r).$

If $X:M \to \r^3$ is a conformal minimal immersion and $\gamma\subset M$ is an oriented closed curve,  the flux  of $X$ on $\gamma$ is given by $p_X(\gamma):=\int_\gamma \mu(s) ds,$ where   $s$ is the oriented arclength parameter on $\gamma$  and $\mu(s)$ is the conormal vector of $X$ at $\gamma(s)$ for all $s.$ Recall that $\mu(s)$ is the unique unit  tangent vector of $X$ at $\gamma(s)$ such that $\{d X(\gamma'(s)),\mu(s)\}$ is a positive basis. Since $X$ is a harmonic map, $p_X(\gamma)$ depends only on the homology class of $\gamma$ and the well defined flux map $p_X:H_1(M,\z) \to \r^3$ is a group homomorphism. 

Our Fusion Theorem asserts the following (see Theorem \ref{th:sequence}):

\begin{quote} 
{\bf Theorem III (Fusion)} {\em Let $M_1,M_2,\ldots$ be a finite or infinite sequence of pairwise disjoint Riemann surfaces with finite conformal type and non empty boundary. For each $n \in \n$ let  $X_n:M_n \to \r^3$ be a conformal complete minimal surface of FTC.

Then, for any $\epsilon >0$   there exist an open parabolic extension  $M^*$ of $\cup_{n} M_n$ and a conformal complete minimal immersion $Y:M^* \to \r^3$ of WFTC  such that $\| X_n-Y|_{M_n}\|_0 \leq \epsilon/n$ and $p_Y|_{{\cal H}_1(M_n,\z)}=p_{X_n},$   where $\|\cdot \|_0$ is the norm of the supremum  on $M_n,$  $n\in \n.$}
\end{quote}
The main tool for proving this theorem has been the algebraic bridge principle given in \cite{approx}. This bridge principle allows good control over the conformal structure, asymptotic behavior and flux map of the resulting surface, and supplies a natural connected sum construction for complete minimal surfaces of finite total curvature in $\r^3.$  Interesting results of this kind can be found in Kapouleas \cite{kapou} and  Yang \cite{yang} works. Theorem III is the core of our existence result for space-filling minimal surfaces. The existence of universal minimal surfaces follows from the separability of the moduli space of complete minimal surfaces with FTC and the Fusion Theorem as well.
%%%%%%%%%%%%%%%%%%%%%%%%
%%%%%%%%%%%%%%%%%%%%%%%
%%%%%%%%%%%%%%%%%%%%%
%%%%%%%%%%%%%%%%%%%%%%

\section{Preliminaries on Riemann surfaces} \label{sec:pre}

As usual, we denote by $\c,$ $\overline{\c}=\c \cup \{\infty\}$ and $\d$ the complex plane, the extended complex plane and the conformal unit disc.

Given a Riemann surface $M,$  $\partial(M)$ will denote the one dimensional topological manifold determined by the boundary points of $M.$ Given $S \subset M,$ write $S^\circ$ and $\overline{S}$ for the topological interior and  the topological closure of $S$  in $M,$ respectively.  A  proper connected subset $S \subset M$ is said to be a {\em region} if it is a topological surface with the induced topology. Open connected subsets of $M-\partial(M)$ are said to be {\em domains} of $M.$ 

Given a point $P\in M-\partial(M)$, we denote by $\mu_P$ the {\em harmonic  measure of $M$ with respect to $P$}. For any  Borel measurable set $I \subset \partial(M),$ $\mu_P(I)=u_I(P),$ where $u_I$ is the unique harmonic function on $M$ vanishing on the ideal boundary of $M$ and satisfying $u_I|_I=1,$ $u_I|_{\partial(M)-I}=0.$ 
\begin{definition}
A Riemann surface $M$ with $\partial(M) \neq \emptyset$ is said to be {\it parabolic} if there exists $P\in M-\partial(M)$ such that $\mu_P$ is full, i.e. $\mu_P(\partial M)=1$. 
If $N$ is an open Riemann surface, $N$ is parabolic in the classical sense if and only if $N-D^\circ$ is parabolic as a Riemann surface with boundary for some (hence for any) closed disc $D \subset N.$
\end{definition}
The fact that $\mu_P$ is full does not depend on the interior point $P;$ this follows from the maximum principle. If $P \in \Omega-\partial(\Omega)$ where $\Omega \subset M$ is a proper region,   we denote by $\mu_P^\Omega$ the harmonic measure relative to $\Omega$ with respect to $P.$ Notice that parabolic surfaces are exactly those on which the maximum principle holds. See \cite{ahlfors} for a good setting.

Let $M$ and $M^*$ be two Riemann surfaces  possibly with non empty boundary.  The surface $M^*$ is said to be an {\em annular extension} of $M$ if $M^*$ is an extension of $M$ and the connected components of $M^*- M^\circ$ are either simply or doubly connected, that is to say, homeomorphic to either $[0,1]\times \s^1$ or  $\overline{\d}-\{E\}$ for some  $E \in \overline{\d}.$ In this case   $j_*:{\cal H}_1(M,\r) \to {\cal H}_1(M^*,\r)$ is an isomorphism.

\subsection{Approximation results on Riemann surfaces}

Let $N$ be a  Riemann surface with $\partial(N)=\emptyset,$ and let $S\subset N$ denote a  subset different from $N$ and satisfying $\overline{S^\circ}\cap S=S$ (for instance, a finite collection of pairwise disjoint regions in $N$).  A connected component $V $ of $N-S$ is said to be {\em bounded} if $\overline{V}$ is compact.

\begin{definition} 
We denote by $\Fg^N(S)$ the space of  functions $f:S \to \overline{\c}$ being meromorphic on some open neighborhood of $S$ in $N.$  If $S$ is open then $\Fg^N(S)$ coincides with the space of meromorphic functions on $S,$ hence does not depend on $N$ and is simply written $\Fg(S).$ 
We write $\Fg^N_0(S)$ (and $\Fg_0(S)$ when $S$ is open) for the corresponding space of holomorphic functions.  
\end{definition}

All these spaces will be endowed with the $\omega(S)$-topology, namely, the topology of the uniform convergence on $S.$ To be more precise, we shall say that a function $f \in \Fg^N(S)$   can be uniformly approximated on $S$ by functions in $\Fg(N)$ if there exists  $\{f_n\}_{n \in \n} \subset \Fg(N)$ such that $\{|f_n|_S-f|\}_{n \in \n} \to 0$ uniformly on $S.$ In particular all $f_n$ have the same set ${\cal P}_f$ of poles  on $S.$ Likewise we define the uniform approximation  of an $f \in \Fg_0^N(S)$ by functions in $\Fg_0(N).$

A complex $1$-form $\theta$ on $S$ is said to be of type $(1,0)$ if for any conformal chart $(U,z)$ on $N$ we have that $\theta|_{U \cap S}=f(z) dz$ for some  $f:U \cap S \to \overline{\c}.$ For instance, holomorphic and meromorphic $1$-forms on $N$ are of type $(1,0).$ 
\begin{definition}
We denote by $\Wg^N(S)$  the space of  $1$-foms of type $(1,0)$ on $S$  being meromorphic on some open neighborhood of $S$ in $N.$ 
If $S$ is open then $\Wg^N(S)$ coincides with the space of meromorphic 1-forms on $S,$ hence does not depend on $N$ and is simply written $\Wg(S).$ 
We write $\Wg^N_0(S)$ (and $\Wg_0(S)$ when $S$ is open) for the corresponding space of holomorphic 1-forms.  
\end{definition}

A 1-form $\theta \in \Wg^N(S)$ can be uniformly approximated on $S$ by 1-forms in $\Wg(N)$ if there exists  $\{\theta_n\}_{n \in \n}\subset \Wg(N)$ such that $\{\frac{\theta_n-\theta}{dz}|_{S \cap K}\}_{n \in \n} \to 0$ in the $\omega(S\cap K)$-topology  for any closed conformal disc $(K,z)$ on $N$ (we are assuming that  $z:K \to z(K)\subset \c$ extends as a  conformal parameter beyond $K$ in $N$).   In particular all $\theta_n$ have the same set  of poles ${\cal P}_\theta$  on $S.$  As above, we say that $\{\theta_n|_S\}_{n \in \n } \to \theta$ in the $\omega(S)$-topology. Likewise we define the uniform approximation of a  $\theta\in\Wg_0^N(S)$ by 1-forms in $\Wg_0(N).$\\

%%%%%%%%%%%%%%%%
%%%%%%%%%%%%%%%
%%%%%%%%%%%%%%%%%%
%%%%%%%%%%%%%%%%%
%%%%%%%%%%%%%%%%%
%%%%%%%%%%%%%%%%
Let $\div(S)$ denote the free commutative group of divisors of $S$ with multiplicative notation. If $D=\prod_{i=1}^n Q_i^{n_i} \in \div(S),$ where $n_i \in \z-\{0\}$ for all $i,$ the set $\{Q_1,\ldots,Q_n\}$ is said to be the {\em support} of $D.$ Let $\deg:\div(S) \to \z$ be the  group homomorphism given by the degree map $\deg(\prod_{j=1}^t Q_j^{n_j})=\sum_{j=1}^t n_j.$ A divisor $D \in \div(S)$ is said to be {\em integral} if $D=\prod_{i=1}^n Q_i^{n_i}$ and $n_i\geq 0$ for all $i.$ Given $D_1,$ $D_2 \in \div(S),$ $D_1 \geq D_2$ if and only if $D_1 D_2^{-1}$ is  integral. For any $f \in \Fg^N(S)$ we denote by $(f)_0$ and $(f)_\infty$ its associated integral divisors of zeroes and poles in $S,$ respectively, and label $(f)=(f)_0/(f)_\infty$ as the divisor associated to $f$ on $S.$ Likewise we define $(\theta)_0$ and $(\theta)_\infty,$ and write  $(\theta)=(\theta)_0/(\theta)_\infty$ for the divisor of $\theta$ in $S,$   $\theta \in \Wg^N(S).$

We need the following extension of  Runge's theorem (see \cite{bishop},  \cite{sche1} and \cite{sche2} for a good setting):

\begin{theorem} \label{th:runge} Let $N$ be a Riemann surface with $\partial(N)=\emptyset,$  and let $S \subset N$ be a finite collection of pairwise disjoint compact regions in $N.$ Let $E \subset N-S$ be a subset meeting each bounded component of $N-S$ in a unique point. 

Then any function $f \in \Fg^N(S)$ can be uniformly approximated on $S$ by functions $\{f_n\}_{n \in \n}$ in $\Fg(N) \cap {\Fg_0}(N-(E\cup {\cal P}_f)),$ where ${\cal P}_f=f^{-1}(\infty)\subset S.$
\end{theorem}

\subsection{Compact Riemann surfaces}\label{subsec:riemann}
The background of the following results can be found, for instance, in \cite{farkas}. 

In the sequel, $R$ will denote an elliptic   Riemann surface of genus $\nu \geq 1.$ 

 Label $H_1(R,\z)$ as the $1^{st}$  homology group with integer coefficients of $R.$ Let $B=\{a_j,b_j\}_{j=1,\ldots,\nu}$ be a canonical homology basis of $H_1(R,\z),$ and write $\{\xi_j\}_{j=1,\ldots,\nu}$ the associated dual basis of $\Wg_0(R),$ that is to say, the one satisfying that $\int_{a_k} \xi_j=\delta_{jk}, \quad j,\;k=1,\ldots, \nu.$ 
 
 Denote by  $\Pi=(\pi_{jk})_{j,\,k=1,\ldots,\nu}$ the Jacobi period matrix with entries $\pi_{jk}=\int_{b_k} \xi_j, \quad j,\;k=1,\ldots, \nu.$ This matrix is symmetric and has positive definite imaginary part. We denote by $L(R)$ the lattice over $\z$ generated by the $2\nu$-columns of the $\nu \times 2 \nu$ matrix $(I_\nu,\Pi),$ where $I_\nu$ is the identity matrix of dimension $\nu.$ 
 
 Finally, call $J(R)=\c^\nu/L(R)$ the Jacobian variety of $R,$ which is a compact, commutative,  complex, $\nu$-dimensional Lie group. 
Given $E \in R,$ denote by $\varphi_{E}: \div(R) \to J(R), \quad \varphi_{E}(\prod_{j=1}^s Q_j^{n_j})= \sum_{j=1}^s n_j \,{}^t(\int_{E}^{Q_j} \xi_1,\ldots,\int_{E}^{Q_j} \xi_\nu)$ the Abel-Jacobi map with base point $E,$ where ${}^t(\, \cdot \,)$ means matrix transpose. If there is no room for ambiguity we simply write $\varphi.$ 

Abel's theorem asserts that  $D\in \div(R)$ is the principal divisor associated to a meromorphic function $f \in \Wg(R)$ if and only if $\deg(D)=0$ and $\varphi(D)=0.$  Jacobi's theorem says that $\varphi:R_\nu \to J(R)$ is surjective and has maximal rank (hence a local biholomorphism) almost everywhere, where $R_\nu$ denotes the space of integral divisors in $\div(R)$ of degree $\nu.$   

Riemann-Roch theorem says that $r(D^{-1})=\deg(D)-g+1+i(D) $ for any $D \in \div(R),$  where $r(D^{-1})$ (respectively, $i(D)$) is the dimension of the complex vectorial space  of functions  $f\in \Fg(R)$ (respectively, 1-forms $\theta \in \Wg(R)$) satisfying that $(f) \geq D^{-1}$ (respectively, $(\theta) \geq D$). 

A point $Q \in R$ is said to be a {\em Weierstrass point} if there exists a non constant meromorphic function $h \in \Fg(R)$ satisfying that $(h)_\infty \leq Q^\nu.$  The number of Weierstrass points in $R$ lies in between $2\nu-2$ and $\nu(\nu^2-1).$

\subsection{Bridge constructions for Riemann surfaces}

A Riemann surface $M$ (possibly with non empty compact boundary) is of finite conformal type if and only if it has finite topology and is parabolic. 
The {\em Osserman compactification} $M^c$ of $M$ is obtained by filling out the conformal punctures corresponding to the topological  ends of $M.$    Moreover, if we attach conformal discs on the holes of $M^c$  we get an elliptic Riemann surface $R$ that we will call a  {\em conformal compactification} of $M.$  Notice that $M^c$ is unique up to biholomorphisms,  whereas $R$ depends on the gluing process.
With this language, $$M^c=R-(\cup_{j=1}^b U_j)\;\;\mbox{and}\;\; M=M^c-\{E_1,\ldots,E_a\},$$ where $U_j,$ $j=1,\ldots,b$ are open discs in $R$ with pairwise disjoint closures in $R$ and $\{E_1,\ldots,E_a\}\subset M^c-\partial(M^c).$

The following notion of conformal connected sum captures some natural bridge constructions for Riemann surfaces. We include the details just for completeness.

Let $M_1,$ $M_2$ be two disjoint Riemann surfaces of finite conformal type and non empty boundary, and fix disjoint Jordan arcs $\gamma_i \subset \partial({M}_i),$ $i=1,2.$ Without loss of generality  assume that ${M}^c_1 \cap {M}^c_2=\emptyset$ as well and write $M=M_1 \cup M_2$ and ${M}^c= {M}^c_1\cup {M}^c_2.$ Let $S$ be a closed conformal disc disjoint from $M^c,$ and introduce a mark  on $S$ consisting of two  distinct Jordan arcs $\gamma_1',$ $\gamma_2' \subset \partial (S).$   By definition,  $\Upsilon=(\{\gamma_1,\gamma_2\},\{S,\gamma_1',\gamma_2'\})$ is said to be a {\em conformal bridge} between $M_1$ and $M_2.$ 

\begin{figure}[h]
\begin{center}
\includegraphics[width=6.9cm,height=4.2cm]{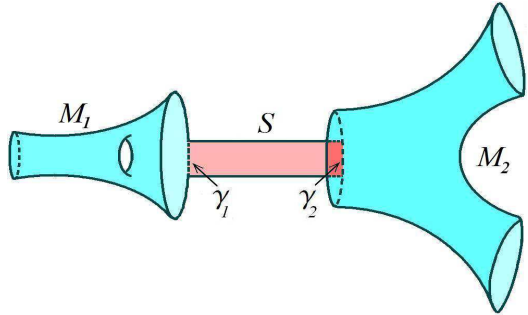} \caption{The surface $M_1\sharp_\Upsilon M_2$}  \label{fig:prime} 
\end{center}
\end{figure}

The surfaces $M_1$ and $M_2$ can be connected via $\Upsilon$  as follows. Take a biholomorphism $w:S \to [0,1] \times [-\delta,\delta]\subset \c$ such that $\gamma'_i=w^{-1}(s_i),$ where  $s_i$ is the segment ${\{i-1}\} \times [-\delta,\delta],$ $i=1,2$ (the real number $\delta$ is uniquely determined by the mark $\{\gamma_1', \gamma_2'\}$ on $\partial(S)$). 
Take a closed disc $V_i \subset M_i$ such that $\gamma_i = V_i  \cap \partial(M_i)$, $i=1,2,$ and  $V_1\cap V_2=\emptyset,$ and consider biholomorphisms   $w_1:V_1 \to [-1,0]\times [-\delta,\delta]$ and $w_2:V_2 \to [1,2]\times [-\delta,\delta]$ such that  $w_i(\gamma_i)=s_i,$ $i=1,2.$ Then simply attach $S$ to $M$ by identifying the points $w^{-1}((i-1,t))$ and $w_i^{-1}((i-1,t))$ for any $t \in [-\delta,\delta] ,$ $i=1,2.$
\begin{definition}
We write $M_1\sharp_\Upsilon M_2$ for the quotient surface,  and call it a {\em connected conformal sum} (or simply, {\em a conformal sum}) of $M_1$ and $M_2$ via $\Upsilon$ (see Figure \ref{fig:prime}). 
\end{definition} 

Up to the projection map to the quotient, $\gamma_i =\gamma_i',$ $i=1,2.$ Furthermore,   $M_1,$ $M_2$  and $S$ are subsets of $M_1\sharp_\Upsilon M_2$ satisfying $M_i \cap S=\gamma_i,$ $i=1,2,$   and $M_1\sharp_\Upsilon M_2= M \cup S.$ Adding the natural chart from $V_1\cup S \cup V_2$ onto $[-1,2]\times [-\delta,\delta]$ induced by $w,$ $w_1$ and $w_2,$  $M_1\sharp_\Upsilon M_2$ becomes a Riemann surface of finite conformal type and non empty boundary and $H_1(M_1\sharp_\Upsilon M_2,\z)=H_1(M_1,\z) \oplus {\cal H}_2(M_2,\z).$  A conformal  compactification  $R_\Upsilon$ of $M_1 \sharp_\Upsilon M_2$ is  said to be {\em a conformal   compactification  of} $M$ via  $\Upsilon.$ Obviously these constructions  guarantee the uniqueness of neither $M_1 \sharp_\Upsilon M_2$ nor $R_\Upsilon,$ because they depend on  the gluing processes.\\

This bridge construction can be used for generating parabolic Riemann surfaces of arbitrary topology. Indeed, let $\{M_j\}_{1\leq j <\sigma},$  $\sigma \in \n \cup \{+\infty\},$ be  sequence  of pairwise disjoint Riemann surfaces of finite conformal type and non empty boundary, and call $M:=\cup_{1\leq j<\sigma} M_j.$ Set $W_1=M_1,$ and working inductively,  for each $j<\sigma$  choose a bridge $\Upsilon_j$ between $W_{j}$ and $M_{j+1}$  and set $W_{j+1}=W_{j} \sharp_{\Upsilon_j} M_{j+1}.$  By definition, the Riemann surface $\sharp_{\Upsilon_j\in {\bf \Upsilon}}M_j: =\cup_{0\leq j<\sigma} W_{j+1}$ is said to be a {\em conformal sum} of $\{M_j\}_{1\leq j <\sigma}$  via the multi-bridge ${\bf \Upsilon}=\{\Upsilon_j\}_{1\leq j<\sigma}.$  Notice that $\sharp_{\Upsilon_j\in {\bf \Upsilon}}M_j$ has genus $\sum_{j<\sigma} \nu_j,$ where $\nu_j$ is the genus of $M_j$ for all $j.$

An open Riemann surface $M^*$ is said to be a {\em parabolic completion} of $M$  via ${\bf \Upsilon}$  if $M^*$ is  parabolic and there exists a proper {\em topological} embedding ${\cal I}: \sharp_{\Upsilon_j\in {\bf  \Upsilon}}M_j \to   M^*$ such that ${\cal I}|_{M}=\mbox{Id}_{M}$ and $M^*$ is an annular extension of  ${\cal I}(\sharp_{\Upsilon_j\in {\bf  \Upsilon}}M_j).$ In particular $M^*$ has genus $\sum_{j<\sigma} \nu_j$ as well and  ${\cal I}_*:H_1(\sharp_{\Upsilon_j\in {\bf  \Upsilon}}M_j,\z) \to H_1(M^*,\z)$ is an isomorphism.  Up to the group monomorphism induced by the inclusion map, $H_1(M_j,\z)$ is a subset of $H_1(M^*,\z),$ $1\leq j<\sigma,$ and therefore $H_1(M^*,\z)$ is the  direct sum $\oplus_{1\leq j<\sigma}H_1(M_j,\z).$ 

\begin{lemma} \label{lem:parabocomple}
 Given $M=\cup_{1\leq j<\sigma} M_j$ and ${\bf \Upsilon}=\{\Upsilon_j\}_{1\leq j<\sigma-1}$ as above, $M$  admits a parabolic completion $M^*$ via ${\bf \Upsilon}.$
\end{lemma}
\begin{proof} Assume first that $\sigma<+\infty.$ In this case $\sharp_{\Upsilon_j\in {\bf  \Upsilon}}M_j$ is of finite conformal type. Let $R$ be the conformal compactification   of $\sharp_{\Upsilon_j\in {\bf  \Upsilon}}M_j$ and consider  a finite subset $E\subset R$ containing all the ends of $M$ and meeting each component of $R-\sharp_{\Upsilon_j\in {\bf  \Upsilon}}M_j$ in a  unique point. It suffices to take $M^*=R-E$ and set ${\cal I}$ as the inclusion map.

Suppose now that $\sigma=+\infty.$ 
Fix a closed disc $D \subset M_1-\partial(M_1)$ and a point $P \in D-\partial(D).$ As above, put $W_1=M_1$ and  $W_{j+1}=W_{j} \sharp_{\Upsilon_{j}} M_{j+1}$ for each $j\geq 1.$ 
Let $c_1^j \subset \partial(W_j)$ and $c_2^j \subset \partial (M_{j+1})$ be the two boundary components (closed Jordan curves) connected by $\Upsilon_j,$ $j\geq 1.$

Let us construct a  sequence $N_1 \subset N_2 \subset \ldots$ of Riemann surfaces  and proper embeddings ${\cal I}_j:W_j \to N_j,$ $j\geq 1,$ such that:
\begin{enumerate}[(a)]
\item  ${\cal I}_{j}|_{W_{j-1}}={\cal I}_{j-1},$ $(\cup_{h > j} M_j ) \cap N_j=\emptyset,$  $\cup_{h \leq j} M_j  \subset N_j$ and  ${\cal I}_j|_{\cup_{h \leq j} M_j}$ is the inclusion map, $j \geq 2,$ 
\item  $N_j$ is a Riemann surface  of   finite conformal type and $\partial(N_{j})$ is a Jordan curve homologically equivalent in $N_j$ to ${\cal I}_j(c_1^{j}),$  $j \geq 1,$  
\item $N_{j}$ is an annular extension of both ${\cal I}_j(W_{j})$ and $N_{j-1} \sharp_{\Upsilon'_{j-1}} M_{j}$ for a suitable bridge $\Upsilon_{j-1}'$ in $N_j$ connecting $\partial(N_{j-1})$ and $c_2^{j-1},$   $j\geq 2.$ 
\item $\mu_P^{{N}_j-D^\circ}(\partial(D))>\frac{j-1}{j},$ where   $\mu_P^{{N}_j-D^\circ}$ is the harmonic measure of  ${N}_j-D^\circ$ with respect to $P,$ $j\geq 1.$
\end{enumerate}
Let $R_1$ be an open parabolic annular extension of  $M_1,$ and notice that $R_1$ is biholomorphic to a finitely punctured compact Riemann surface.  Without loss of generality suppose that $R_1 \cap (\cup_{h>1} M_h)=\emptyset.$ Since $R_1$ is parabolic,  we can find a  proper  region $N_1 \subset R_1$ such that $N_1$ has just one hole (hence connected boundary), $M_1 \subset N_1-\partial(N_1),$ $\partial(N_1)$ is homologically equivalent $c_1^1,$ and $N_1$ is an annular extension of  $M_1.$ Set ${\cal I}_1:M_1 \to N_1$ the inclusion map and observe that $\mu_P^{{N}_1-D^\circ}(\partial(D))>0.$ The above items hold for $j=1.$ 

Reasoning inductively, suppose that  we have constructed ${N}_j$ and ${\cal I}_j,$ $1\leq j \leq m-1,$ satisfying the above properties. Take a bridge  $\Upsilon'_{m-1}$ between  $N_{m-1}$ and $M_m$ connecting $\partial(N_{m-1})$ and $c_2^{m-1}.$ Let $R_m$ be an open parabolic annular extension of $N_{m-1} \sharp_{\Upsilon'_{m-1}} M_m.$ As above notice that $R_m$ is a finitely punctured compact Riemann surface, and without loss of generality suppose that $R_m \cap (\cup_{h>m} M_h) =\emptyset.$  Let ${\cal I}_m:W_m \to R_m$ be any extension of ${\cal I}_{m-1}$ as a proper topological embedding satisfying that  ${\cal I}_m|_{M_m}=\mbox{Id}_{M_m}.$  Since $R_m$ is parabolic, there exists a proper   region $N_m \subset R_m$ with just one hole  such that ${\cal I}_m(W_m)\subset N_m-\partial(N_m),$  $N_m$ is an annular  extension of both ${\cal I}_m(W_m)$ and $N_{m-1} \sharp_{\Upsilon'_{m-1}} M_m,$   $\partial(N_m)$ is homologically equivalent to ${\cal I}_m(c_1^m)$ and  $\mu_p^{{N}_m-D^\circ}(\partial(D))>1-1/m.$ Considering the natural embedding  ${\cal I}_m:W_m \to N_m,$  the induction is closed.

Set $M^*=\cup_{j\geq 1} N_j$ and ${\cal I}:\sharp_{\Upsilon_j\in {\bf  \Upsilon}} M_j\to M^*,$ ${\cal I}|_{W_j}={\cal I}_j$ for all $j.$ It is not hard to check that the open Riemann surface $M^*$ is an annular extension of ${\cal I}(\sharp_{\Upsilon_j\in {\bf  \Upsilon}}M_j).$ Moreover, if $\mu_P^{M^*-D^\circ}$ is the harmonic measure of  $M^*-D^\circ$ with respect to $P$ then $\mu_p^{M^*-D^\circ}(\partial(D))=\lim_{j \to \infty} \mu_p^{N_j-D^\circ}(\partial(D))=1,$ proving that  $M^*$ is parabolic and the lemma.
\end{proof}

\section{Preliminaires on minimal surfaces} \label{sec:appro}

Throoughout this section, $N$ will be an open  Riemann surface and $M \subset N$ a  finite union of pairwise disjoint regions with compact boundary.

\begin{definition}\label{def:basica} 
\begin{enumerate}[(a)]
Let $ \mathcal{E}(N)$ denote the space of conformal complete minimal immersions $X:N \to\r^3$ of WFTC. Likewise, we write $\mathcal{E}_N(M)$ for the space of conformal complete minimal immersions $X:M \to\r^3$ of WFTC  that extend as a conformal minimal immersion to some neighborhood of $M$ in $N.$
  \end{enumerate}
\end{definition}
If $M$ is open then ${\cal E}(M)={\cal E}_M(M).$ When $M$ has finite conformal type,  $ \mathcal{E}_N(M)$ is the space of  conformal complete minimal immersions of $M$ in  $\r^3$ with  FTC that  extend to some neighborhood of $M$ in $N.$ These spaces will be endowed with the following ${\cal C}^0$ topology:

\begin{definition} \label{def:convergencia}
A sequence $\{X_n\}_{n \in \n} \subset  \mathcal{E}_N(M)$ is said to converge in the ${\cal C}^0(M)$-topology to  $X_0 \in \mathcal{E}_N(M)$ if for any proper region $\Omega \subset M$ of finite conformal type, $\{X_n|_\Omega\}_{n \in \n} \to X_0|_\Omega$ uniformly on $\Omega.$  If $M$ has finite conformal type, this topology coincides with the one of uniform convergence on $M.$
\end{definition}

Take $X\in \mathcal{E}_N(M)$ and write $(\phi_1,\phi_2,\phi_3)$ for the complex differential  $\partial_z X.$ Notice that $\partial_z X\in \Wg^N_0(M)^3.$ Since $X$ is conformal and minimal,  then  $\phi_1=\frac{1}{2}(1/g-g) \phi_3$ and $\phi_2=\frac{i}{2}(1/g+g) \phi_3,$ where $g \in \Fg^N(M)$ and up to the stereographic projection coincides with the Gauss map of $X.$ The pair $(g,\phi_3)$ is known as the {\em Weierstrass representation} of $X$ (see \cite{osserman}).

Clearly $X(P)=X(Q)+\mbox{Re} \int_{Q}^P (\phi_1,\phi_2,\phi_3),$ $P,$ $Q \in M.$ The induced intrinsic metric $ds^2$ on $M$ and its Gauss  curvature $\mathcal{K}$ are given by the expressions: 
\begin{equation} \label{eq:regular}
ds^2=\sum_{j=1}^3 |\phi_j|^3=\frac{1}{4} |\phi_3|^2(\frac{1}{|g|}+|g|)^2,\quad \mathcal{K}=-\left(\frac{4 |dg||g|}{ |\phi_3| (1+|g|^2)^2}\right)^2.
\end{equation}
The total curvature of $X$ is given by $c(X):=\int_M \mathcal{K} dA,$ where  $dA$ is the area element of $ds^2,$ and the {\em flux map} of $X$ by  the expression $p_X:H_1(M,\z) \to \r^3,$ $p_X(\gamma)=\mbox{Im} \int_\gamma \partial_z X.$

By Huber, Osserman and Jorge-Meeks results \cite{huber,osserman,jorge-meeks}, if $X$ is complete and of FTC then $X$ is proper,  $M$ has finite conformal type and the Weierstrass data of $X$ extend meromorphically to $M^c.$
 
\begin{remark} Take $\{X_n,\; n \in \n\} \cup \{X\} \subset \mathcal{E}_N(M)$ and assume that  $\{X_n\}_{n \in \n}\to X$ in $\mathcal{E}_N(M).$ If $\Omega \subset M$ is a proper region of finite conformal type,  it is not hard to see that $\{(X_n- X)|_{{\Omega}^c}\}_{n \in \n} \to 0$ uniformly on  ${\Omega}^c$ and the Weierstrass data of $X_n$ converge in the $\omega({\Omega}^c)$-topology to the ones of $X$ on ${\Omega}^c.$ Indeed, just observe that by Riemann's removable singularity theorem,  $X_n-X$ extends harmonically to the punctures of $\Omega$ for all $n$ and  $\{X_n-X\}_{n \in \n} \to 0$ uniformly on ${\Omega}^c$ as well.
\end{remark}

\begin{figure}[h]
\begin{center}
\includegraphics[width=5cm,height=2.5cm]{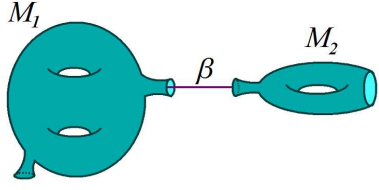} \caption{The surfaces $M_1,$ $M_2$ and the curve $\beta.$}  \label{fig:segun} 
\end{center}
\end{figure}

\begin{figure}[h]
\begin{center}
\includegraphics[width=7.5cm,height=6.3cm]{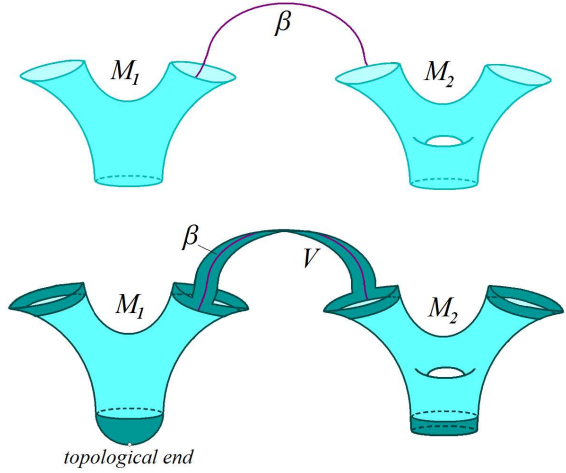} \caption{An annular neighborhood $V$ of $M\cup \beta.$}  \label{fig:terce} 
\end{center}
\end{figure}

Let  $M_1,$ $M_2$ be two proper regions in $N$  with finite conformal type  and non-empty boundary, call $M=M_1 \cup M_2$ and consider  an analytic regular Jordan arc $\beta \subset N$ with endpoints $P_1\in \partial(M_1)$ and $P_2\in \partial(M_2)$ and  otherwise disjoint from $M$ (see Figure \ref{fig:segun}). We will always assume that $\beta$ lies in an  open analytic arc $\beta_0$ in $N$ such that  $\beta_0-\beta \subset M-\partial(M).$ 
A proper region $V \subset N$ is said to be an {\em annular extension} of $M \cup \beta$ in $N$  if $M \cup \beta\subset V^\circ,$  $V-(M^\circ\cup \beta)$ contains no connected components with compact closure that are disjoint from $\partial(V),$  and $V-(M^\circ\cup \beta)$ consists of a finite collection of compact annulus and once punctured closed discs. In particular, the induced homomorphism  $j_*:H_1(M \cup \beta,\z) \to H_1(V,\z)$ is an isomorphism, where $j:M\cup \beta \to V$ is the inclusion map. See Figure \ref{fig:terce}.
A map $X:M\cup \beta \to \r^3$ is said to be {\em smooth} is $X|_{M_j},$ $j=1,2,$  and $X|_{\beta_0}$ are smooth.
For instance, if $Y:N \to \r^3$ is a smooth map then $X=Y|_{M \cup \beta}$ is smooth.
\begin{definition} We denote by ${\cal E}_N(M \cup \beta)$ the space of smooth maps $X:M\cup \beta \to \r^3$ such that $X_j:=X|_{M_j}\in \mathcal{E}_N(M_j),$ $j=1,2,$  and $X|_{\beta}$ is a smooth immersion. This space is endowed with the ${\cal C}^0(M \cup \beta)$-topology of the uniform convergence on $M \cup \beta.$
\end{definition} 
It is clear that $Y|_{M \cup \beta}\in {\cal E}_N(M \cup \beta)$  for all $Y \in \mathcal{E}(N).$

Notice that $H_1(M \cup \beta,\z)= H_1(M,\z)=H_1(M_1,\z) \oplus H_1(M_2,\z),$ and for each $X\in {\cal E}(M\cup \beta)$ identify the flux map  $p_{X}:H_1(M \cup \beta,\z) \to \r^3$ of $X$   with the one of $X|_M.$

The proof of the following results can be found in \cite{approx}:
%%%%%%%%%%%%%%%%%%%%%%%%%%%%%%%
%%%%%%%%%%%%%%%%%%%%%%%%%%%%%
\begin{theorem}[The Algebraic Bridge Principle] \label{th:density} Assume that $N$ is an open Riemann surface of finite conformal type and $N-(M \cup \beta)$ consists of a finite collection of pairwise disjoint once punctured conformal discs. Let $X$ be an arbitrary immersion in $\mathcal{E}_{N}(M\cup \beta).$

Then there exists a sequence  $\{Y_n\}_{n \in \n}  \subset \mathcal{E}(N)$  such that $p_{Y_n}|_{H_1(M \cup \beta)}=p_{X}$ for all $n \in \n$ and  $\{Y_n|_{M\cup \beta}\}_{n \in \n} \to X$ in the ${\cal C}^0(M \cup \beta)$-topology. 
Furthermore, if $C$ is a positive constant and  $V$ an annular extension of $M\cup \beta$ in $N,$   then  $\{Y_n\}_{n \in \n}$ can be chosen in such a way that $d_{Y_n}(M\cup \beta, \partial(V)) \geq C$ for all $n,$ where $d_{Y_n}$ is the intrinsic distance in $N$ induced by $Y_n.$
\end{theorem}

\begin{theorem}[General Approximation]  \label{th:parabo}
Let $M$ be a Riemann surface of finite conformal type and  $\partial(M)\neq \emptyset,$  and  let $M^*$ be an extension of $M$ with $\partial(M^*)=\emptyset.$   Consider an immersion $X\in {\cal E}_{M^*}(M)$ and a linear extension $q:H_1(M^*,\z) \to \r$ of  $p_X.$ 

Then, there exists a  sequence of conformal complete minimal immersions $\{Y_n\}_{n \in \n}\in {\cal E}(M^*)$ such that  $\{Y_n|_M\}_{n \in \n} \to Y$ in the ${\cal C}^0(M)$-topology and $p_{Y_n}=q.$
\end{theorem}

\section{Fusion theorems for minimal surfaces of FTC}\label{sec:fusion}
Parabolicity is a powerful tool because it ensures the well-posedness (existence, uniqueness, stability...) of interesting geometrical problems.  Theorem \ref{th:parabo} implies the existence of complete minimal surfaces of WFTC with arbitrarily prescribed (non compact) parabolic conformal structure. In particular,  there exist complete parabolic minimal surfaces with arbitrary topology. Our interest resides in obtaining fusion theorems for this kind of surfaces.

Given a topological space $T$ and a continuous map $f:T \to \r^3,$ we write  $\|f\|_0=\sup \{\|f(P)\| \,:\, P \in T\},$ where $\|\cdot\|$ is the Euclidean norm.

\begin{theorem}[Fusion] \label{th:sequence}
Let $\{M_j\}_{1\leq j<\sigma}$ be a  sequence of pairwise disjoint Riemann surfaces of finite conformal type and {\em non empty}   boundary, where $\sigma \in \n \cup \{+\infty\},$ and let  $M^*$ be a parabolic completion of $M=\cup_{1\leq j<\sigma}M_j.$ 
Consider  $X_i \in {\cal E}_{M^*}(M_i),$ $i \geq 1,$ and fix $\epsilon>0.$

Then there is  $Y \in {\cal E}(M^*)$  such that  $\|Y|_{M_j}-X_j\|_0 \leq \epsilon/j$  and $p_Y|_{H_1(M_j,\z)}=p_{X_j}$ for all $j\geq 1.$
\end{theorem} 
\begin{proof} Consider the multi-bridge ${\bf \Upsilon}=\{\Upsilon_j\}_{1\leq j<\sigma-1}$ and  proper embedding ${\cal I}: \sharp_{\Upsilon_j\in {\bf  \Upsilon}}M_j \to   M^*$ such that ${\cal I}|_{M}=\mbox{Id}_{M}$ and $M^*$ is an annular extension of  ${\cal I}(\sharp_{\Upsilon_j\in {\bf \Upsilon}}M_j).$ 
Up to natural identifications, we will assume that ${\cal I}$ is the inclusion map and $M^*$ is an annular extension of a conformal sum $\sharp_{\Upsilon_j\in {\bf  \Upsilon}} M_j.$ Fix a point $P \in M_1.$

 Like in the proof of Lemma \ref{lem:parabocomple}, we can find an exhaustion $N_1 \subset N_2 \subset \ldots$ of $M^*$ by  proper regions of of finite conformal type and (compact)  connected boundary, and  a sequence of bridges $\Upsilon_1',$ $\Upsilon_2',\ldots$ between $N_1$ and $M_2,$ $N_2$ and $M_3,\ldots,$ 
 such that:  $N_1$ is an annular extension of $M_1,$ $\cup_{h\leq j} M_j \subset  N_j,$ $(\cup_{h> j} M_j)\cap  N_j=\emptyset$   and $N_{j+1}$ is an annular extension of $N_{j}\sharp_{\Upsilon'_{j}}M_{j+1},$ $j\geq 1.$   When $\sigma<+\infty$ the sequence ends at $N_{\sigma-1}=M^*.$

Let us construct  $Y_j \in {\cal E}_{M^*}(N_j),$ $1\leq j<\sigma,$ such that:    $\|Y_j|_{M_{j}}-X_{j}\|_0\leq \frac{\epsilon}{j 2^j},$ $p_{Y_j}|_{H_1(M_{j},\z)}=p_{X_{j}}$   and $d\big(Y_j(P), Y_{j} \big(\partial({N}_j)\big)\big) >j$ for all $j\geq 1,$ and $\|Y_j|_{N_{j-1}}-Y_{j-1}\|_0 \leq \frac{\epsilon}{j 2^j}$  and $p_{Y_j}|_{H_1(N_{j-1},\z)}=p_{Y_{j-1}}$  for all $j\geq 2.$    

Indeed, by Theorem  \ref{th:parabo}, there is $Y_1\in {\cal E}_{M^*}(N_1)$ such that $\|Y_1|_{M_1}-X_1\|_0\leq \epsilon/2,$ $p_{Y_1}|_{H_1(M_1,\z)}=p_{X_1}$ and $d\big(Y_1(P), Y_1 \big(\partial({N}_1)\big)\big) >1.$ Reasoning inductively, suppose we have constructed $Y_j\in {\cal E}_{M^*}(N_j)$ satisfying the above properties.  By Theorem \ref{th:density}, there exist  $Y_{j+1}\in {\cal E}_{M^*}(N_{j+1})$ such that   $\|Y_{j+1}|_{M_{j+1}}-X_{j+1}\|_0,$ $\|Y_{j+1}|_{N_{j}}-Y_{j}\|_0\leq \epsilon/(j+1) 2^{j+1},$  $p_{Y_{j+1}}|_{H_1(M_{j+1},\z)}=p_{X_{j+1}},$ $p_{Y_{j+1}}|_{H_1(N_{j},\z)}=p_{Y_{j}}$ and   $d\big(Y_{j+1}(P), Y_{j+1} \big(\partial({N}_{j+1})\big)\big) >j+1,$ closing the induction. 

When $\sigma<+\infty$ the immersion  $Y=Y_{\sigma-1} \in {\cal E}(M^*)$ solves the theorem.

If $\sigma=+\infty,$ there exists a possibly branched conformal  minimal immersion $Y:M^*\to\r^3$ such that $\{Y_m|_{N_j}\}_{m \in \n} \to  Y|_{N_j}$ uniformly on $N_j$ for any $j\geq 1.$ Furthermore,  $\|Y|_{M_j}- X_j\|_0 \leq \epsilon/j$  and $\|Y|_{N_j}- Y_j\|_0 \leq \epsilon/(j+1)$ for any $j\geq 1.$ Let us show that $Y$ has no branch points. Indeed, let $(g_j,\phi_3^j)$ denote the Weierstrass data of $Y_j,$ $j \geq 1,$ and likewise call $(g,\phi_3)$ as the ones of $Y.$ Obviously, $\{g_j,\phi_3^j)\}_{j \in \n} \to (g,\phi_3)$ uniformly on compact subsets of $M^*.$  Take an arbitrary $P_0 \in M^*,$ and consider $j_0 \in \n$ such that  $P_0 \in N_{j_0}^\circ.$ Up to a rigid motion, $g(P_0) \neq 0,$ $\infty,$ hence we can find an closed disc $D \subset N_{j_0}$ such that $P_0 \in D^\circ$ and the functions $g|_D$ and  $g_j|_D,$ $j \in \n,$ are holomorphic and without zeroes. Since $Y_j$ has no branch points, $\phi_3^j$ has no zeroes on $D$ for all $j.$  By Hurwith theorem, either $\phi_3=0$ of $\phi_3$ has no zeroes on $D$ as well. In the first case the identity principle gives $\phi_3=0$ on $M^*,$ contradicting that $\|Y-X_1\|_0<\epsilon$ on $M_1$ when $\epsilon$ is small enough, and proving that  $Y_\epsilon$ is an immersion.

Finally, let us see that $Y$ is complete and of WFTC. By Osserman's theorem, the Gauss map of $Y_j$ extends meromorphically to  ${N}^c_j,$  $j \in \n.$  Since $\|Y_j-Y|_{{N}_j}\|_0$ is finite then $Y_j-Y$ extends harmonically to ${N}^c_j$ and  $Y|_{N_j}$ is complete and of finite total curvature  for any $j.$  It remains to check that $Y$ is complete. First notice that those curves in $M^*$ diverging to an annular  end of some  $N_j$ have infinite intrinsic length with respect to $Y.$ Moreover, the fact that  $d\big(Y_j(P), Y_{j} \big(\partial({N}_j)\big)\big) >j$ for all $j$ implies that  $\lim_{j \in \n} d\big(Y(P), Y \big(\partial({N}_j)\big)\big)=+\infty.$ This shows that any  divergent curve in $M^*$ that does  not diverge to an annular end of some $N_j$ have infinite intrinsic length as well, and so $Y \in \mathcal{E}(M^*).$ Clearly   $p_Y|_{H_1(M_j,\z)}=p_{X_j}$ for all $j\geq 1$ and we are done. 
\end{proof}

This fusion theorem can be used for producing minimal surfaces with exotic geometry. We start with the following existence result for space-filling minimal surfaces:

\begin{corollary} \label{co:space-filling}
For each $\nu \in \{0\}\cup \n \cup \{\infty\},$ there exists a space-filling, open,  parabolic, complete and minimal surface in $\r^3$ with genus $\nu,$  WFTC and no symmetries. 
\end{corollary}
\begin{proof} Let $\{v_1,v_2,v_3\}\subset \r^3$ be three linearly independent unit  vectors in general position, that is to say, such that $\langle v_{i_1},v_{j_1} \rangle \neq \pm \langle v_{i_2},v_{j_2} \rangle$ provided that $\{i_1,j_1\} \neq \{i_2,j_2\},$ where $\langle,\rangle$ is the Euclidean metric. Let $\{r_n:\, n \in \n\}$ be a bijective enumeration of $\q,$ and write $\Sigma_{i,n}=r_n v_i+\{u \in \r^3 :\, \langle u,v_i\rangle=0\;\mbox{and}\; \langle u,u\rangle \geq 1/n^2\},$  $i=1,2,3,$ $n \in \n.$ Consider a conformal parameterization $Y_{i,n}:M_{i,n} \to \r^3$ of $\Sigma_{i,n},$ where $M_{i,n} \cap M_{j,m}=\emptyset$ provided that $(i,n) \neq (j,m).$ 

Let $Z:N \to \r^3$ denote a conformal parameterization of the Chen-Gackstatter genus one minimal surface, take a closed disc $D \subset N$ and write $Z_0=Z|_{N_0},$ where $N_0=N-D^\circ.$ Recall that $N$ (hence $N_0$) has an only topological end, and $Z_0(N_0)$ is asymptotic to the classical Enneper surface. In particular, $Z_0(N_0)$ is not asymptotic to a plane. Call $N_{0,n}=N_0 \times \{n\}$ and set $Z_{0,n}:N_{0,n} \to \r^3,$ $Z_{0,n}((P,n))=Z_0(P)$ for all $n<\nu+1.$  Let $\{Y_j:N_j \to \r^3 :\; j \in \n\}$ denote a bijective enumeration of $\{Y_{i,n}\;:\; i=1,2,3,\, n \in \n\} \cup \{Z_{0,n} \;:\; n < \nu+1\}.$ 

Let $N_j^*$ be an annular neighborhood of $N_j$ homeomorphic to $N_j,$   and without loss of generality suppose that $Y_j$ can be extended to $N_j^*,$ $j \in \n.$  
 Take  a parabolic completion $M^*$ of $\{N_j^*\}_{j \in \n}$ and observe that $Y_j \in {\cal E}_{M^*}(N_j)$ for all $j.$ Consider the fusion immersion $Y\in {\cal E}( M^*)$ associated to  $\{Y_j\}_{j \in \n}$ and any $\epsilon>0$ via  Theorem \ref{th:sequence}. 

It is not hard to check that $Y$ is space-filling (we leave the details to the reader). Moreover $M^*$ has genus $\nu,$ hence it suffices to check that $Y$ has no symmetries. Reason by contradiction, and suppose there exists a rigid motion $\sigma:\r^3 \to \r^3$ different from the indentity map $\mbox{Id}$ and leaving invariant $Y(M^*).$ Call $\sigma_0:M^* \to M^*$ as the  intrinsic isometry satisfying that $Y \circ \sigma_0=\sigma \circ Y.$ The embedded planar annular ends of $Y$ have limit normal vector parallel to some $v_i,$ $i \in \{1,2,3\}.$ As $\sigma_0$ maps annular ends onto annular ends with the same geometry, then $\vec{\sigma}$ leaves invariant the system of vectors $\{\pm v_j\,:\; j=1,2,3\},$ where $\vec{\sigma}$ is the linear isometry associated to $\sigma.$ Taking into account that the vectors $v_1,$ $v_2$ and $v_3$ are placed in general position, we infer that  $\vec{\sigma}=\pm \mbox{Id}.$

Assume for a moment that $\vec{\sigma}=\mbox{Id},$ that is to say,  $\sigma$ is a non trivial translation. In this case $\sigma_0$ takes annular planar ends on annular planar ends with the same limit normal vector. Fix $i \in \{1,2,3\},$ and for each $n \in \n$ let $m(n) \in \n$ denote the unique natural number such that $M_{i,m(n)}$  and $\sigma_0(M_{i,n})$ determine the same annular end. Call  $\Omega_{i,n}:=\sigma_0(M_{i,n})\cap M_{i,m(n)}$ and notice that the Euclidean distance $d(Y(\Omega_{i,n})- r_{m(n)} v_i,\Sigma_i) \to 0$ as $n \to \infty,$ where $\Sigma_i=\{u  \;:\; \langle u,v_i\rangle=0,\; u \neq 0\}.$ This clearly implies that the translation vector of $\sigma$ must be  orthogonal to $v_i.$  However, this can not occur for all $i\in \{1,2,3\},$ getting a contradiction.

Suppose now that $\vec{\sigma}=-\mbox{Id},$ i.e.,  $\sigma$ is a symmetry with respect to a point $P_0 \in \r^3.$ Reasoning as above, $\sigma$ preserves the annular ends of Enneper type. However, the Enneper type ends of $Y(M^*)$ lie in a neighborhood of radius $\epsilon$ of the genus one Chen-Gakstatter surface $Z(N),$  and this surface has no central symmetries. This contradiction concludes the proof.
\end{proof} 
\subsection{Universal minimal surfaces}\label{subsec:univ}
This section is devoted to the existence problem of universal minimal surfaces. 

We start with some notations. 
Let $M$ be a Riemann surface of finite conformal type with $\partial(M) \neq \emptyset,$ and let $N$ be an open Riemann surface.  An immersion $Y\in {\cal E}(N)$  is said to {\em pass by} $X\in {\cal E}(M)$ if there exist  proper  regions $\{\Omega_n\}_{n \in \n}$  in $N$ and biholomorphisms $h_n:M \to \Omega_n,$ $n \in \n,$   such that  $\{Y \circ h_n\}_{n \in \n} \to X$ in the ${\cal C}^0(M)$-topology. Note that if $Y$  passes by $X$ then  $X(M)\subset \overline{Y(N)},$ but the converse is not necessarily true.

\begin{definition}
Let $N$ be an open Riemann surface $N.$ An immersion $Y\in {\cal E}(N)$ is said to be {\em universal} if for any compact Riemann surface $M$ with non empty boundary and any conformal minimal immersion $X:M \to \r^3,$   $Y$  passes by $X.$ 
\end{definition}
The next lemma is an elementary consequence of Theorem \ref{th:sequence}:

\begin{lemma} \label{lem:sequence}
Let $\{Y_i:N_i \to \r^3\}_{1\leq i<\sigma}$  be a sequence of conformal complete minimal immersions of FTC, where $\sigma\in \n \cup \{+\infty\},$ and assume that $\partial(N_i)\neq \emptyset$ for all $i.$  

Then there exists  an open  parabolic Riemann surface $M^*$ and an immersion $Y \in {\cal E}(M^*)$ such that $Y$ passes by $Y_i$ for all $i.$
\end{lemma}
\begin{proof} Recall that $N_i$ has finite conformal type for all $i,$ and without loss of generality assume that $N_i \cap N_j=\emptyset,$ $i \neq j.$  Set $N_{i,j}=N_i \times\{j\},$ $h_{i,j}: N_i \to N_{i,j},$ $h_{i,j}(P)=(P,j),$ $Y_{i,j}=Y_i \circ h_{i,j}^{-1}$ for all $j \in \n.$ 
Write ${\cal Y}$ for the countable family $\{Y_{i,j}:\; 1\leq i<\sigma,\; j \in \n\},$ and take  a bijective enumeration $\{X_j:M_j \to \r^3\}_{j \in \n}$ of ${\cal Y}.$ Label $M_j^*$ as an annular neighborhood of $M_j$ homeomorphic to $M_j^*$  where $X_j$ can be extended as a conformal minimal immersion, $j \in \n,$ and let $M^*$ denote a parabolic completion of $\{M_j^*\}_{j \in \n}.$  Note that $X_j \in {\cal E}_{M^*}(M_j)$ for all $j,$ and consider the fusion immersion  $Y \in {\cal E}(M^*)$ of Theorem \ref{th:sequence} associated to $\{X_i\}_{i \in \n}$ and $\epsilon>0.$ For any $i<\sigma$ and $j \in \n,$ label $i_j$ as the unique natural such that $Y_{i,j}=X_{i_j}$ (hence $N_{i,j}=M_{i_j}$). As $\lim_{j \to \infty} \|Y|_{M_{i_j}}-X_{i_j}\|_0=0,$ then $ \{Y\circ h_{i_j}\}_{j \in \n} \to Y_i$ in the ${\cal C}^0(N_i)$-topology, concluding the proof. 
\end{proof}

In order to approach the existence problem of  universal minimal surfaces, we need some preliminary results on Riemann surfaces. We start with the following:

\begin{lemma} \label{lem:branching}
Let $R$ be an elliptic Riemann surface, and let $V$ be an open disc in $R.$ Then  there is an $f_V \in \Fg(R)$ all of whose branch points lie in $V.$
\end{lemma}
\begin{proof} The proof is trivial when $R=\overline{\c}.$ Then we will assume that $R$ has positive genus $\nu.$ 
For the following, it is convenient to go over again the notations and results of Section \ref{sec:pre}, and specially those of Subsection \ref{subsec:riemann}.

\begin{assertion}\label{ass:prime}
There exists  a $\tau_0\in \Wg(R) \cap \Wg^R_0(R-V)$ without zeroes in $R-V.$
\end{assertion}
\begin{proof}
Fix  $E \in V$ and  take a non zero  $\theta\in \Wg_0(R).$ Put $(\theta)=D_1 \cdot D,$  where $D_1 \in \div(R-V)$ and $D \in \div(V).$   By Jacobi's theorem, we can find an open disc $U \subset V$  such that $\varphi_E:U_\nu \to \varphi_E(U_\nu)$ is a diffeomorphism, where $U_\nu$ is the set of divisors in $R_\nu$  with support in $U.$ Since $J(R)$ is a compact additive Lie Group and $\varphi_E(U_\nu) \subset J(R)$ is an open subset,  one has $n_0 \varphi_E(U_\nu)=J(R)$ for large enough $n_0 \in \n$ Therefore, there is $D_2 \in U_\nu$ such that $\varphi_E(D_2^{n_0})=\varphi_E(D_1)=\varphi_E(D_1 E^m),$ where $m=n_0 \nu -\deg(D_1).$ By  Abel's theorem there exists $f_0 \in \Fg(R)$ such that $(f_0)=\frac{D_2^{n_0}}{D_1 E^\nu}.$ It suffices to set $\tau_0= f_0 \,\theta.$ 
\end{proof}

Fix a non Weierstrass point $Q \in V,$ and label ${\cal U}_Q \subset \Wg(R)$ as the complex vectorial subspace of meromorphic 1-forms with $(\theta)\geq Q^{-\nu-1}.$ By Riemann-Roch theorem, $dim_{\c} \;{\cal U}_Q=2 \nu$ and  the map ${\cal G}:{\cal U}_Q \to \c^{2 \nu},$ ${\cal G}(\tau)=(\int_{c} \tau)_{c \in B},$ is a linear isomorphism.

As usual write $B=\{a_j,b_j\}_{1\leq 1\leq \nu}$ for a canonical basis of $H_1(R,\z),$ and choose the representative curves $a_j,$ $b_j,$ $j=1,\ldots,\nu,$ in $R-\overline{V}.$

\begin{assertion} \label{ass:log}  Let $W\subset R$ be an open disc containing $\overline{V}$ and disjoint from $a_j,$ $b_j$ for all $j.$ Then, for any function $h\in \Fg^R_0(R-W)$ never vanishing on $R-W,$ there exists $f \in \Fg_0(R-\{Q\})$  never vanishing on $R-\{Q\}$   such that $\log(h/f)$ has a well defined branch on $R-W.$
\end{assertion}
\begin{proof} Take  $\tau \in {\cal U}_Q$ such that $dh/h-\tau$ has vanishing periods along $a_j,$ $b_j$ for all $j,$ and observe that $\frac{1}{2 \pi i}\int_{a_j} \tau,$ $\frac{1}{2 \pi i}\int_{b_j} \tau \in \z$ for all $j.$ Set $h_0=\int (dh/h-\tau)\in \Fg^R_0(R-W)$ and $f=e^{\int\tau}\in \Fg_0(R-\{Q\}).$ Finally, note that $f$ never vanishes on $R-\{Q\}$ and  $\log(h/f)=h_0\in \Fg^R_0(R-W).$
\end{proof}

Let $\sigma$ be a non null  exact 1-form in $\Wg_0(R-\{Q\}),$ and  let $W\subset R$ be an open disc  containing  $\overline{V}$ and all the zeroes of $\sigma_0$ in $R-\{Q\}.$
\begin{assertion} \label{ass:segun}
There exists $\kappa \in \Wg(R-\{Q\}) \cap \Wg^R_0(R-V)$ without zeroes on $R-V$  and $g_0\in {\cal F}_0^R(R-W)$ such that $\sigma|_{R-W}= e^{g_0} (\kappa|_{R-W}).$
\end{assertion}
\begin{proof}
Set $h=(\sigma/\tau_0)|_{R-W},$ where $\tau_0$ is the 1-form  given in Claim \ref{ass:prime}. If necessary, choose the representative curves $a_j,$ $b_j,$ $j=1,\ldots,\nu,$ for $B$ in $R-W.$ By the previous claim, there is $f \in \Fg_0(R-\{Q\})$ never vanishing on $R-\{Q\}$   such that $g_0:=\log(h/f)$ is a well defined holomorphic map on $R-W.$ Label $\kappa=f \tau_0 \in \Wg(R-\{Q\}) \cap \Wg^R_0(R-V),$ and note that $\kappa$ has no zeroes on $R-V.$ Finally, observe that $\sigma|_{R-W}= (e^{g_0} \kappa)|_{R-W}\in \Wg^R_0(R-W).$
\end{proof}

\begin{assertion} \label{ass:sobre}
The linear map ${\cal L}_0:\Fg_0 (R-\{Q\}) \to \c^{2 \nu},$ ${\cal L}_0(h)= (\int_c h e^{g_0} \kappa)_{c \in B},$ is surjective.
\end{assertion}
\begin{proof} Endow $\Fg_0 (R-\{Q\})$ with the topology of the uniform convergence {\em on compact subsets} of $R-\{Q\},$ and observe that ${\cal L}_0$ is continuous. Take a basis $\{\theta_j\}_{j=1,\ldots,{2 \nu}}$ of ${\cal U}_Q,$  set $h_j=\theta_j/((e^{g_0} \kappa) \in \Fg^R_0(R-V)$ for each $j,$ and observe that  $\{(\int_c h_j e^{g_0} \kappa)_{c \in B}\}_{j=1,\ldots,2 \nu}$ is a basis of $\c^n.$

On the other hand, Theorem \ref{th:runge} implies that $\Fg_0(R-\{Q\})$ is dense in $\Fg^R_0(R-W)$ with respect to the $\omega(R-W)$-topology, and so  
$h_j$ lies in the closure of $\Fg_0(R-\{Q\})$  in $\Fg^R_0(R-W)$ for all $j.$ By a continuity argument ${\cal L}_0$ is surjective and we are done.
\end{proof}
Consider $\{g_n\}_{n \in \n} \subset \Fg_0(R-\{Q\})$ such that $\{g_n|_{R-W}\}_{n \in \n}\to g_0|_{R-W}$ in the $\omega(R-W)$-topology (use Theorem \ref{th:runge}). Set  ${\cal L}_n:\Fg_0 (R-\{Q\}) \to \c^{2 \nu},$ ${\cal L}(h)= (\int_c h e^{g_n} \kappa)_{c \in B},$  and observe that ${\cal L}_n$ is a continuous linear operator for all $n \in \n \cup \{0\}.$ Furthermore, $\{{\cal L}_n\}_{n \in \n} \to {\cal L}_0$ in the weak topology, that is to say, $\{{\cal L}_n(h)\}_{n \in \n}\to {\cal L}_0(h)$ for all $h \in \Fg_0 (R-\{Q\}).$  By Claim \ref{ass:sobre} there exists $\{f_j\}_{j=1,\ldots,2 \nu}\subset  \Fg_0(R-\{Q\})$ such that $\{{\cal L}_n(f_j)\}_{j=1,\ldots,2 \nu}$ generates $\c^{2 \nu},$ $n$ large enough (up to removing finitely many terms, for all $n\in \n\cup\{0\}$).  
Define ${\cal Q}_n:\c^{2 \nu} \to \c^{2 \nu},\quad {\cal Q}_n(\{x_j\}_{j=1,\ldots,2 \nu}) = (\int_c  e^{g_n+\sum_{j=1}^{2 \nu} x_j f_j} \kappa)_{c \in B},$ $n \in \n \cup \{0\},$ and notice that $\{{\cal Q}_n\}_{n \in \n} \to {\cal Q}_0$ as analytic maps on compact subsets of $\c^{2\nu}.$ Since the Jacobian of ${\cal Q}_n$ at ${\bf 0}=(0)_{j=1,\ldots,2\nu}$ is different from zero  for all $n \in \n \cup \{0\},$   there exists an Euclidean ball $B_0 \subset \c^{2 \nu}$ centered at ${\bf 0}$ such that ${\cal Q}_n:B_0 \to {\cal Q}_n(B_0)$ is a diffeomorphism for all $n \in \n\cup \{0\}.$ Furthermore, as ${\cal Q}_0({\bf 0})={\bf 0}$ then ${\bf 0} \in {\cal Q}_n(B_0)$ for  large enough $n$ (without loss of generality, for all $n$). Take $(y_j^n)_{j=1,\ldots,2\nu}\in B_0$ such that ${\cal Q}_n((y_j^n)_{j=1,\ldots,2\nu})={\bf 0}$ and set $\sigma_n=e^{g_n+\sum_{j=1}^{2 \nu} y_j f_j} \kappa\in \Wg_0(R-\{Q\}),$ $n \in \n.$  The 1-form $\sigma_n$  have no periods and never vanish on  $R-V,$ hence the  function $F_n=\int \sigma_n\in \Fg^R_0(R-V)$ has no branch points on $R-V,$ $n \in \n.$

To finish, fix $n_0 \in \n$ and use Theorem \ref{th:runge} to find $\{H_k \}_{k \in \n} \subset \Fg(R) \cap \Fg_0(R-\{Q\})$ such that $\{H_k \}_{k \in \n}\to F_{n_0}$ in the $\omega(R-V)$-topology. By Hurwitz theorem, we can suppose that $d H_n$ never vanishes on $R-V$ for all $n.$ It suffices to choose $f_V=H_n$ for some $n \in \n.$ 
\end{proof}

Given a polynomial $\pg$ with complex coefficients in the variables $z$ and $w,$ we denote by $\mbox{Deg}_z(\pg)$ and $\mbox{Deg}_w(\pg)$  the degree of  $\pg$ in $z$ and $w,$ respectively.

Let $R$ be an elliptic Riemann surface of genus $\nu\geq 1.$  For any $f \in {\cal F}(R),$ write  $\mbox{Deg}(f)$ for the degree of $f$ as meromorphic function on $R.$
Let $Q \in R$ be a non Weierstrass point, and for each $n \geq \nu+1,$ let $f_n\in \Fg(R) \cap \Fg_0(R-\{Q\})$ denote a non zero function with $Deg(f_n)=n$ and  polar divisor $(f_n)_\infty=Q^n.$  Label $z=f_{\nu+1}$ and $w=f_{\nu+2}.$ We know there is an irreducible complex  polynomial $\pg$ in the variables $z$ and $w$  with $\mbox{Deg}_z(\pg)-1=\mbox{Deg}_w(\pg)=\nu+1$   satisfying $\pg(z(P),w(P))=0$ for all $P \in R.$ Furthermore,  $R$ is biholomorphic to the algebraic curve $C_\pg:=\{(z,w) \in \overline{\c}^2:\, \pg(z,w)=0\}$ (up to this biholomorphism we will consider $R=C_\pg$), and the pair $\{z,w\}$ generates the field of meromorphic functions ${\cal F}(R).$ The last means that any $f \in \Fg(R)$ is of the form $f=\pg_1(z,w)/\pg_2(z,w)$ for suitable polynomials $\pg_1$ and $\pg_2$  without common factors and having $\mbox{Deg}_w(\pg_i) \leq \nu.$ 

\begin{remark}
If $R=\overline{\c},$ we also have that $R \cong C_{\pg_0}:=\{(z,w) \in \overline{\c}^2 :\, \pg_0(z,w)=0\}$ for $\pg_0(z,w)=w^2-(z-a_1)(z-a_2),$ where $a_1,$ $a_2 \in \c,$ $a_1 \neq a_2.$
\end{remark}

\begin{definition} For each $v=(\nu,k,s) \in (\n \cup \{0\}) \times \n^2,$ write ${\cal W}_{v}$ for the space of couples $(\pg,F)$ such that:
\begin{itemize}
\item $\pg(z,w)$ is an irreducible complex polynomial in $(z,w)$ with $\mbox{Deg}_z(\pg)-1=\mbox{Deg}_w(\pg)=\nu+1.$
\item The algebraic curve $C_\pg$ has genus $\nu,$ $(\infty,\infty)$ is the only pole of $z$ and $w$ as meromorphic functions on $C_\pg,$  and $(0,0) \in C_\pg.$ 
\item $F=\big((\pg_{1,j},\pg_{2,j})\big)_{j=1,2,3},$ where $\pg_{1,j}$ and $\pg_{2,j}$ are complex polynomials in $(z,w)$ with no common factors so that $\mbox{Deg}_w(\pg_{i,j}) \leq \nu,$ $i=1,2,$ $j=1,2,3.$ 
\item Setting  $f_j:C_\pg \to \overline{\c},$ $f_j(P):=\pg_{1,j}(z(P),w(P))/\pg_{2,j}(z(P),w(P)),$ we have that $\sum_{j=1}^3 f_j^2=0$ on $C_\pg$  and $\mbox{Deg}(g)=k,$ where $g:=f_3/(f_1-i f_2).$
\item  The polar set $E_{\pg,F}$ of the vectorial 1-form $F dz$ on $C_\pg$ has $s$ points, $(0,0) \notin E_{\pg,F},$ and $\sum_{j=1}^3 |f_j|^2 |dz|^2$ has no zeroes on $C_\pg-E_{\pg,F}.$ 
\item The meromorphic 1-form $f_j dz\in \Wg(R)$ has no real periods on $C_\pg-E_{\pg,F},$ $j=1,2,3.$
\end{itemize}
We also set  ${\cal A}_{v}= \r^3 \times {\cal W}_{v}.$
\end{definition}

For any two complex polynomials $\pg_1(z,w)=\sum_{i,j} a_{i,j}z^i w^j$ and $\pg_2(z,w)=\sum_{i,j} b_{i,j}z^i w^j,$ we  set $d(\pg_1,\pg_2)=\sum_{i,j} |a_{i,j}-b_{i,j}|.$ We endow ${\cal W}_{v}$   with the topology induced by the metric  $d^7\equiv d \times (d\times d)^3,$ and likewise equip ${\cal A}_{v}= \r^3 \times {\cal W}_{v}$ with the topology induced by the metric $d_0\times d^7,$ where $d_0$ is the Euclidean metric in $\r^3.$

Given $v=(\nu,k,s)$ and $(\pg,F) \in {\cal W}_v$ as above, elementary algebraic arguments show that $\mbox{Deg}_z(\pg_{h,j}),$ $h=1,2,$ $j=1,2,3,$ admit an universal upper bound depending only on $k$ and $\nu.$ Notice also that for any $y=(x,(\pg,F)) \in {\cal A}_v,$ the well defined immersion $$X_y:C_p-E_{\pg,F}\to \r^3,\quad  X_y(q)=x+\mbox{Re}\big(\int_{(0,0)}^q  F dz\big)$$ lies in ${\cal E}(C_\pg-E_{\pg,F})$ for all $(\pg,F) \in {\cal W}_v.$

For each  $v=(\nu,k,s) \in (\n \cup \{0\}) \times \n^2,$ let ${\cal E}_{v}$ denote the moduli space of conformal complete minimal immersions $X:M \to \r^3$ such that $M$ is a $s$-punctured genus $\nu$ elliptic Riemann surface  and $X$ has total curvature $-4 \pi k.$ It is clear  that $X_y \in {\cal E}_v$ for any $y \in {\cal A}_v.$ 

Equip  ${\cal E}_{v}$ with the following topology: a sequence $\{X_n:M_n \to \r^3\}_{n \in \n} \subset {\cal E}_v$ is said to be convergent in the ${\cal C}^0_*$ topology to $X_0:M_0 \to \r^3,$ where $X_0 \in {\cal E}_v,$  if for any compact region $\Omega_0 \subset M_0$ there exist compact regions $\Omega_n \subset M_n$ and biholomorphisms $h_n:\Omega_0 \to \Omega_n,$ $n \in \n,$ such that $\{X_n \circ h_n\}_{n \in \n} \to X_0|_{\Omega_0}$ in the ${\cal C}^0(\Omega_0)$-topology.

\begin{lemma} \label{lem:topologia}
The map $\Delta_v:{\cal A}_v \to {\cal E}_v,$ $\Delta_{v}(y):=X_y$  is surjective and continuous.
\end{lemma}
\begin{proof} For the surjectivity, take an arbitrary immersion $X:M \to \r^3$ in ${\cal E}_v.$  By Osserman's theorem,  $M=R-\{Q_1,\ldots,Q_s\},$ where $R$ is an elliptic genus $\nu$ Riemann surface, and the Weierstrass data of $X$ extend meromorphically to $R.$  Fix a non Weierstrass point $Q \in M$ and as above take $z,$ $w \in \Fg(R)$ with $\mbox{Deg}(z)=\mbox{Deg}(w)-1=\nu+1$ and $(z)_\infty=Q^{\nu+1},$ $(w)_\infty=Q^{\nu+2}.$ Fix also $Q_0 \in M-\{Q\},$ and without loss of generality suppose $z(Q_0)=w(Q_0)=0.$ Label  $\pg(z,w)$ as the irreducible polynomial in $(z,w)$ such that $R=C_\pg,$ and write $\partial_z X/dz=(f_j(z,w))_{j=1,2,3},$ where $f_j\in \Fg(R)$ is a rational function of the form $\pg_{1,j}(z,w)/\pg_{2,j}(z,w)$ and  $\mbox{Deg}_w(\pg_{i,j}) \leq \nu,$ $i=1,2,$ $j=1,2,3.$ As the meromorphic Gauss map $g=f_3/(f_1-i f_2)$ has degree $k$ on $R,$ then $\big(\pg,F=(\pg_{1,j},\pg_{2,j})_{j=1,2,3}\big) \in {\cal W}_{v}.$ It is clear that $X=\Delta_{v} \big( (X(Q_0),(\pg,F)) \big).$

To check that  $\Delta_v$ is continuous, take $\{y_n=(x_n,(\pg_n,F_n))\}_{n \in \n\cup\{0\}} \in {\cal A}_v$ such that $\{y_n\}_{n \in \n} \to y_0,$ and fix an arbitrary compact region $\Omega_0 \subset C_{\pg_0}-E_{\pg_0,F_0}.$  We have  to find compact regions $\Omega_n \subset C_{\pg_n}-E_{\pg_n,F_n}$ and biholomorphisms $h_n:\Omega_0 \to \Omega_n,$ such that $\{X_{y_n}\circ h_n\}_{n \in \n} \to X_{y_0}|_{\Omega_0}$ in the ${\cal C}^0(\Omega_0)$-topology.
\begin{remark} Recall that $z$ and $w$ are meromorphic functions on  $C_{\pg_n}=\{(z,w) \in \overline{c}^2: \pg_n(z,w)=0\}$ {\em for all} $n\in \n\cup\{0\},$ hence  they depend on $n.$ This lack of notation does not affect our exposition.   
\end{remark}

Let $U_0$ be an open subset of $C_{\pg_0}-(\{(0,0)\}\cup \Omega_0)$ containing $E_{\pg_0,F_0}.$  By  Lemma \ref{lem:branching}, there is a meromorphic function $z_0:C_{\pg_0} \to \overline{\c}$ with all its branch points in $U_{0}.$ Write $z_0=\qg_1(z,w)/\qg_2(z,w),$ where $\qg_1,$ $\qg_2$ are polynomials with no common factors and $\mbox{Deg}_w(\qg_i) \leq \nu,$ $i=1,2,$ and choose $w_0$ any function in $\{z,w\}$ so that $\{z_0,w_0\}$ generates $\Fg(C_{\pg_0}).$ Let $z^l w^j$ be the effective monomial (i.e., with non zero coefficient) in $\qg_1$ and $\qg_2$ with maximum degree {\em as meromorphic function on} $C_{\pg_0}.$ Since $\mbox{Deg}(w)-1=\mbox{Deg}(z)=\nu+1$ and   $\mbox{Deg}_w(\qg_i) \leq \nu,$ $i=1,2,$ this monomial always exists and is unique. Furthermore, as $z$ and $w$ have an unique pole at the same point (namely, $(\infty,\infty)$) of $C_{\pg_0},$  then $\mbox{Deg}(z_0)=l (\nu+1)+ j (\nu+2).$ The same argument shows that $z_n:C_{\pg_n}\to \overline{\c},$   $z_n=\qg_1(z,w)/\qg_2(z,w),$ has  $\mbox{Deg}(z_n)=\mbox{Deg}(z_0)$ as meromorphic function on  $C_{\pg_n}$ for $n$ large enough (without loss of generality, for all  $n \in \n$). 

In the sequel we write  $a=\mbox{Deg}(z_n)$ (which does not depend on $n$) and $E_n=E_{\pg_n,F_n},$  $n \in \{0\} \cup \n.$  We also label $B_n$ as the branch point set of $z_n$ on $C_{\pg_n}$ for all $n \in \{0\} \cup \n.$  For any $P \in C_{\pg_0},$ denote  $b_P$ as the branching number of $z_0:C_{\pg_0} \to \overline{\c}$ at $P,$ and for each $\zeta \in \overline{\c}$ write $a_{\zeta}=\sum_{P \in z_0^{-1}(\zeta)} b_{P}.$ 

Choose for each $\zeta \in z_0(B_0\cup E_0)$ an open disc $D_\zeta \subset \overline{\c}$ centered at $\zeta$ so that:
\begin{itemize}
\item   $\{\overline{D}_\zeta:\,\zeta \in z_0(B_0\cup E_0)\}$ is a family of pairwise disjoint closed discs, 
\item $z_0^{-1}(D_\zeta)$ consists of $a-a_\zeta$ conformal discs, 
\item  if $P \in z_0^{-1}(\zeta)$ and  $U_P$ is the connected component of $z_0^{-1}(D_\zeta)$ containing $P,$ then $z_0|_{U_P}:U_P \to D_\zeta$ is a branched covering of $b_P$ sheets, and
\item $\overline{U}_P \subset U_{0}$  when $P \in E_0\cup B_0.$ 
\end{itemize}

Since $\{d(\pg_n,\pg_0)\}_{n \in \n} \to 0,$  $B_0$ is the limit set of $B_n$ as $n \to \infty$  in $\overline{\c}^2.$ In other words, if $\{P_n\}_{n \in \n} \subset \overline{\c}^2$ is a convergent sequence such that  $P_n \in B_n$ for all $n,$ then its limit lies in $B_0,$ and any point of $B_0$ is the limit of a sequence of this kind. Likewise, if we write $F_n=\big((\pg^n_{1,j},\pg^n_{2,j})\big)_{j=1,2,3}$ one has $\{d(\pg^n_{i,j},\pg^0_{i,j})\}_{n \in \n} \to 0$ for all $i,$ $j,$ and so $E_{0}$ is the limit set of $\{E_{n}\}_{n \in \n}$  in $\overline{\c}^2$ as well. 

For each $P \in z_0^{-1}(z_0(B_0 \cup E_0)),$   choose $Q_P^n \in z_n^{-1}(z_0(P)),$ $n \in \n ,$  so that $\{Q_P^n\}_{n \in \n} \to P$ as points of $\overline{\c}^2.$ By elementary topology, and up to removing finitely many terms of the sequence $\{y_n\}_{n \in \n}$ if necessary, we can suppose that:
\begin{enumerate}[(i)]
\item  $z_n^{-1}(D_\zeta)$ is a collection of $a-a_\zeta$ pairwise disjoint open discs on $C_{\pg_n}$ for all $\zeta \in z_0(B_0\cup E_0)$ and $n \in \n\cup\{0\}.$
\item For any $P \in z_0^{-1}(z_0(B_0 \cup E_0))$ and $n \in \n\cup \{0\},$  $z_n|_{U^n_P}:U^n_P \to D_{z_0(P)}$ is a branched covering of $b_P$ sheets, where $U^n_P$ is the component of $z_n^{-1}(D_{z_0(P)})$ containing $Q_P^n.$ \footnote{Notice that $U_P=U_P^0$ for all $P \in z_0^{-1}(z_0(B_0 \cup E_0)).$}
\end{enumerate}

Set $W=\overline{\c}-\cup_{\zeta \in z_0(B_0\cup E_0)} D_\zeta,$  $W_n=z_n^{-1}(W)\subset C_{\pg_n}$ and  $\pi_n:=z_n|_{W_n}:W_n \to W,$ $n \in \n \cup \{0\}.$ 
Fix $\zeta_0 \in W$ and choose $P_0^n \in z_n^{-1}(\zeta_0),$ $n \in \n \cup \{0\},$  so that $\{P_0^n\}_{n \in \n} \to P_0^0$ as points of $\overline{\c}^2.$ Basic monodromy arguments give that $(\pi_n)_*(\Pi_1(W_n))=(\pi_0)_* (\Pi_1(W_0))\subset \Pi_1(W),$ where 
$\Pi_1(W_n)$ is the fundamental group of $W_n$ with base point $P_0^n,$ $\Pi_1(W)$ is the one of $W$ with base point $\zeta_0,$ and $(\pi_n)_*:\Pi_1(W_n) \to \Pi_1(W)$ is the  group homomorphism induced by $\pi_n,$ $n \in \n.$

For each $n \in \n,$ let $\lambda_n:W_0 \to W_n$ denote the unique biholomorphism such that $\lambda_n(P_0^n)=P_0^0$ and $z_n \circ \lambda_n=z_0|_{W_0}.$

Label $J=\{P \in z_0^{-1}(B_0 \cup E_0): b_P=0\},$ and notice that $z_n|_{U^n_P}:U^n_P \to D_{z_0(P)}$ is a biholomorphism for all $P \in J$ and $n \in \n.$  Call $V_n=\cup_{P  \in J} U^n_P$ and $\hat{W_n}=W_n \cup V_n,$  $n \in \n \cup \{0\}.$  Let $\hat{\lambda}_n: \hat{W}_0 \to \hat{W}_n$ denote the natural  extension of $\lambda_n$ satisfying $z_n \circ \hat{\lambda}_n=z_0|_{\hat{W}_0},$ $n \in \n.$ 
Since $\cup_{P \in E_0\cup B_0}\overline{U}_P \subset U_{0},$ then $\Omega_0 \subset \hat{W}_0-\partial(\hat{W}_0).$ Set $\Omega_n=\hat{\lambda}_n(\Omega_0)$ and write $h_n=\hat{\lambda}_n|_{\Omega_0}:\Omega_0 \to \Omega_n,$ $n \in \n.$  
The facts that $\{d(\pg_n,\pg_0)\}_{n \in \n} \to 0$  and  $\{d(\pg^n_{i,j},\pg^0_{i,j})\}_{n \in \n} \to 0$ for all $i,$ $j,$ imply that $\{w \circ \hat{\lambda}_n\}_{n\in \n} \to w|_{\hat{W}_0}$ and  $\{F_n(z\circ \hat{\lambda}_n,w \circ \hat{\lambda}_n)\}_{n \in \n} \to F_0(z|_{\hat{W}_0},w|_{\hat{W}_0})$ in the $\omega(\hat{W}_0)$-topology.  Taking into account that $\{x_n\}_{n \in \n} \to x_0,$  we  deduce that $\{X_{y_n} \circ h_n\}_{n \in \n} \to X_{y_0}|_{\Omega_0}$ in the ${\cal C}^0(\Omega_0)$-topology,  concluding the proof.
\end{proof}

Now we can state the main result of this subsection.

\begin{theorem}[Existence of universal surfaces]  \label{th:universal}
There exist parabolic complete universal minimal surfaces of WFTC in $\r^3.$
\end{theorem}
\begin{proof}Set ${\cal A}=\cup_{v \in (\n \cup \{0\})\times \n^2} {\cal A}_{v}$ and ${\cal E}=\cup_{v \in (\n \cup \{0\})\times \n^2} {\cal E}_v$ endowed with the corresponding direct sum topologies,  and define $\Delta:{\cal A} \to {\cal E},\quad \Delta|_{{\cal A}_{v}}=\Delta_{v}.$

Notice that ${\cal A}_{v}$ is separable,  take  a dense countable subset ${\cal D}_{v} \subset {\cal A}_{v}$ and denote by ${\cal S}_{v}=\Delta_{v} ({\cal D}_{v}).$ Lemma \ref{lem:topologia} says that  ${\cal S}:=\cup_{v \in (\n \cup \{0\})\times \n\times \n}  {\cal S}_{v}$ is a dense countable subset of  ${\cal E}$ as well. 

For each $X:M \to \r^3$ in ${\cal S},$ fix a countable basis $B_X$ of the topology on $M$  formed by open discs bounded by Jordan curves in $M,$  and call ${\cal S}_X=\{X|_{M-D}\,: D \in B_X\}.$ Finally set ${\cal S}_0=\cup_{X \in {\cal S}} {\cal S}_X.$ 

By Lemma \ref{lem:sequence}, there exists an open parabolic Riemann surface $M^*$  and an immersion $Y\in {\cal E}(M^*)$ passing by $X$ for all $X \in {\cal S}_0.$ Let us show that $Y$ is universal.

Let $M_0$ be a compact genus $\nu$ Riemann surface with non empty  boundary, label $s>0$ as the number of components in $\partial(M_0).$ Let $X_0:M_0 \to \r^3$ be a conformal minimal immersion that extends as a conformal minimal immersion to some open Riemann surface $N$ containing $M_0.$ Let $M_0^*$ be a compact annular extension of $M_0$ in $N,$ and construct a  conformal compactification $R$ of $M_0^*.$ Consider a finite subset $E\subset R-M_0$ so that $R-E$ is an annular extension of $M_0$ and notice that $X_0\in {\cal E}_{R-E}(M_0).$ Then take a sequence $\{X_n\}_{n \in \n}\subset {\cal E}(R-E)$ converging to $X_0$ in the  ${\cal C}^0(M_0)$-topology (use Theorem \ref{th:parabo}). Note that $X_n \in {\cal E}_{v_n},$ where  $v_n=(\nu,k_n,s)$ for some $k_n\in \n.$ 
Fix $Q_0 \in M_0,$ and use Lemma \ref{lem:topologia} to find $y_n=(X_n(Q_0),(\pg_n,F_n)) \in {\cal A}_{v_n}$ such that $X_n=\Delta_{v_n}(y_n).$    By the density of ${\cal S}_{v_n}$ in ${\cal E}_{v_n}$ (see Lemma \ref{lem:topologia}), there exists $\{\hat{X}_{j,n}:N_{j,n} \to \r^3\}_{j \in \n} \subset {\cal S}_{v_n},$ regions $W_{j,n} \subset N_{j,n}$ and biholomorphisms $h_{j,n}:M_0 \to W_{j,n},$ $j \in \n,$ such that $\{\hat{X}_{j,n}\circ h_{j,n}\}_{j \in \n} \to X_n|_{M_0}$ in the $C^0(M_0)$-topology. Choose  a disc $D_{j,n} \in B_{X_{j,n}}$ disjoint from $W_{j,n},$ call $X_{j,n}=\hat{X}_{j,n}|_{N_{j,n}-D_{j,n}} \in {\cal S}_{X_{j,n}} \subset {\cal S}_0$ and observe that $\{X_{j,n}\circ h_{j,n}\}_{j \in \n} \to X_n|_{M_0}$ in the $C^0(M_0)$-topology too. Finally, take $j_n\in \n$ such that $\|X_{j_n,n} \circ h_{j_n,n}-X_n|_{M_0}\|_0<1/n$ and label  $h_{n}=h_{j_n,n},$ $n \in \n.$ Since  $\{X_{j_n,n} \circ h_n\}_{n \in \n} \to X_0$ in the ${\cal C}^0(M_0)$-topology and  $Y$ passes by $X_{j_n,n}$ for all $n,$ then $Y$ passes by $X_0$  and we are done. \end{proof}

{\bf FRANCISCO J. LOPEZ} \newline
Departamento de Geometr\'{\i}a y Topolog\'{\i}a \newline
Facultad de Ciencias, Universidad de Granada \newline
18071 - GRANADA (SPAIN) \newline
e-mail: fjlopez@ugr.es

\end{document}